\documentclass[11pt,a4paper,reqno]{amsart}
\usepackage{amsfonts}
\usepackage{graphicx}   
\usepackage{graphics}
\usepackage{epsfig}
\input{psfig.sty}
\usepackage{anysize}
\marginsize{2.5cm}{2.5cm}{2.5cm}{2.5cm}

\numberwithin{equation}{section}

\def\({\left(}
\def\){\right)}

\begin{document}
\setcounter{page}{1}
\bigskip
\bigskip
\title[]
 {The characterization problem for one class of second order
 operator pencil with complex periodic coefficients.}
\author[]{Efendiev R.F. \dag }
\thanks{\dag Institute Applied Mathematics, Baku State University
Z.Khalilov, 23, AZ1148, Baku, Azerbaijan. \underline
rakibaz@yahoo.com }

\begin{abstract}
\bigskip
The purpose of the present work is solving the characterization
problem, which consists of identification of necessary and
sufficient conditions on the scattering data ensuring that the
reconstructed potential belongs to a particular class

 \noindent  MSC: 34B25, 34L05, 34L25,47A40,
81U40

\end{abstract}
\maketitle

\section {Introduction}

The purpose of the present work is solving the characterization
problem, which consists of identification of necessary and
sufficient conditions on the scattering data ensuring that the
reconstructed potential belongs to a particular class. In our case
the potential belongs to $Q_{ +} ^{2} $ consisting of functions

\begin{equation}
\label{eq1}
q\left( {x} \right) = \sum\limits_{n = 1}^{\infty}  {q_{n} exp\left( {inx}
\right)} \quad ,
\end{equation}

\noindent which is a subclass of the class $Q^{2}$of all $2\pi $
periodic complex valued functions on the real axis R, belonging to
$ L_2 [0,2\pi ] $.

The object under consideration is the operator $L$ given by the differential
expression

\begin{equation}
\label{eq2}
l\left( {\frac{{d}}{{dx}},\lambda}  \right) \equiv -
\frac{{d^{2}}}{{dx^{2}}} + 2\lambda p\left( {x} \right) + q\left( {x}
\right) - \lambda ^{2}
\end{equation}

\noindent in the $L_{2} \left( { - \infty ,\infty}  \right)$ with
potentials $p\left( {x} \right) = \sum\limits_{n = 1}^{\infty}
{p_{n} e^{inx}} ,q\left( {x} \right) = \sum\limits_{n =
1}^{\infty}  {q_{n} e^{inx}} $, for which fullfield $
\sum\limits_{n = 1}^\infty  {n \cdot \left| {p_n } \right| <
\infty ,}  \sum\limits_{n = 1}^{\infty} {\left| {q_{n}} \right| <
\infty} , \quad \lambda $-a spectral parameter.

The inverse problem for the potentials (\ref{eq1}) was formulated
and solved in [1, 2], where is shown, that the equation $Ly = 0$,
has the solution

\begin{equation}
\label{eq3}
e_{ \pm}  \left( {x,\lambda}  \right) = e^{ \pm i\lambda x}\left( {1 +
\sum\limits_{n = 1}^{\infty}  {V_{n}^{ \pm}  e^{inx} + \sum\limits_{n =
1}^{\infty}  {\sum\limits_{\alpha = n}^{\infty}  {\frac{{V_{n\alpha} ^{ \pm
}} }{{n \pm 2\lambda} }e^{i\alpha x}}} } }  \right) \quad ,
\end{equation}

$ V_n^ \pm  ,V_{n\alpha }^ \pm $ are some numbers and Wronscian of
the system of solutions is equal to $2i\lambda $.

The limit

 $e_{n}^{ \pm}  \left( {x} \right) = \mathop {lim}\limits_{\lambda \to \mp
n/2} \left( {n \pm 2\lambda}  \right)e_{ \pm}  \left( {x,\lambda}
\right) = \sum\limits_{\alpha = n}^{\infty}  {V_{n\alpha} ^{ \pm}
e^{i\alpha x}e^{ - i\frac{{n}}{{2}}x},\,\,\,\,n \in N} $

\noindent is also a solution of the equation $L{y}= 0$, but is
lineary dependent with $e_{ \pm}  \left( {x, \pm \frac{{n}}{{2}}}
\right)$. Thus, there exist numbers $ \hat S_n ,n \in N, $ for
which

\begin{equation}
\label{eq3}
 e_{n}^{ \pm}  \left( {x} \right) = \hat {S}_{n}^{ \pm}  e_{ \mp}  \left( {x,
\mp \frac{{n}}{{2}}} \right),\,\,\,n \in N.
\end{equation}

 From last relation one
may obtain that$\hat {S}_{n}^{ \pm}  = V_{nn}^{ \pm} $.

In the work [1] the spectral analysis of the operator pencil
$L$was carried out and sufficient condition for reconstruction
$p\left( {x} \right),q\left( {x} \right) \in Q_{ +} ^{2} $ using
the values $ \hat S_n^ \pm  ,n \in N $ was found.

Note that some characterizations for the Sturm-Liouville operator
in the class of real-valued potentials belonging to $L_{1}^{1}
\left( {R} \right)$($L_{\alpha} ^{1} \left( {R} \right)$is the
class of measurable potentials satisfying the condition $
\int\limits_R {dx(1 + \left| x \right|} )^\alpha  \left| {p_\gamma
(x)} \right| < \infty $), have been given by Melin [3] and
Marchenko [4]. (See for details [5, 6, and 7]). For the potential
$ p(x)=0 , q\left( {x} \right) \in Q_{ + }^{2} $ which in the
nontrivial cases is complex valued the inverse problem was first
formulated and solved by Gasymov M.G [8]. Later the complete
solution of the inverse problem for the cases $ p(x) = 0,q(x) \in
Q_ + ^2 $ was found by Pastur and Tkachenko [9].

Let us formulate now the basic result of the present work.

\underline {Definition.} Constructed by the help of formula (1.4)
sequence $\{ \hat {S}_{n}^{ \pm}  \} _{n = 1}^{\infty}  $, is
called as a set of spectral data for the operator $L$  given by
the differential expression (\ref{eq2}) with potential $p\left(
{x} \right),q\left( {x} \right) \in Q_{ +} ^{2} $ .

\underline {Thoerem 1}. In order the given sequence of complex
numbers $\{ \hat {S}_{n}^{ \pm}  \} _{n = 1}^{\infty}  $ to be a
spectral data for the operator $L$  given by the differential
expression (\ref{eq2}) with potentials $p\left( {x}
\right),q\left( {x} \right) \in Q_{ +} ^{2} $ it is necessary and
sufficient, that the conditions are fulfilled at the same time:

1)\begin{equation} \label{eq5}
 \{ n^2 \hat S_n^ \pm  \} _{n = 1}^\infty   \in l_1 ;
\end{equation}

2) Infinite determinant

\begin{equation}
\label{eq4}
D\left( {z} \right) \equiv det\left\| {\delta _{nm} - \sum\limits_{k =
1}^{\infty}  {\frac{{4\hat {S}_{m}^{ -}  \hat {S}_{k}^{ +} } }{{\left( {m +
k} \right)\left( {n + k} \right)}}e^{i\frac{{m + k}}{{2}}z}e^{i\frac{{n +
k}}{{2}}z}}}  \right\|_{n,m = 1}^{\infty}
\end{equation}

\noindent exists, is continuous, is not equal to zero in the
closed semiplane $\overline {C_{ +} }  = \{ z:Imz \ge 0\} $ and is
analytical inside open semiplane $C_{ +}  = \{ z:Imz > 0\} .$
\section { On an inverse problem of the scattering theory in the
semiaxis}

On the base of the proof of the Theorem 2, lies investigation of the
equation$Ly = 0$. Taking in it

\begin{equation}
\label{eq5} x=it,  \lambda=-ik,  y(x)=Y(t)
\end{equation}

\noindent
we get

\begin{equation}
\label{eq6}
 - {Y}''\left( {t} \right) + 2\mu \bar {p}\left( {it} \right)Y\left( {t}
\right) + \bar {q}\left( {it} \right)Y\left( {t} \right) = \mu ^{2}Y\left(
{t} \right)
\end{equation}

\noindent
where

\begin{equation}
\label{eq7} \overline {p} \left( {t} \right) = ip\left( {it}
\right) = i\sum\limits_{n = 1}^{\infty}  {p_{n} e^{ - nt}} ,
\sum\limits_{n = 1}^\infty  {n \cdot \left| {p_n } \right| <
\infty ,} \quad \overline {q} \left( {t} \right) = - q\left( {it}
\right) = - \sum\limits_{n = 1}^{\infty}  {q_{n} e^{ - nt}}
\sum\limits_{n = 1}^\infty  {\left| {q_n } \right| < \infty }  .
\end{equation}

As a result the equation (\ref{eq6}) is obtained, the potential of
which decreases at$t \to \infty $. The specification of the
considered inverse problem is defined by the fact that the
potentials belong to the class$ Q_ + ^2 $. In this section we
suppose $t \in R^{ +} $.

The procedure of analytic continuation that allows from the result for the
(\ref{eq6}) to get corresponding results for the equation (\ref{eq2}) will be
investigated in next section.

The equation (\ref{eq6}) with potentials (\ref{eq7}) has the solution

 \begin{equation}
\label{eq4} f_{ \pm}  \left( {t,\mu}  \right) = e^{ \pm i\mu
t}\left( {1 + \sum\limits_{n = 1}^{\infty}  {V_{n}^{ \pm}  e^{ -
nt} + \sum\limits_{n = 1}^{\infty}  {\sum\limits_{\alpha =
n}^{\infty}  {\frac{{V_{n\alpha} ^{ \pm }} }{{in \pm 2\mu} }e^{ -
\alpha t}}} } } \right)\end{equation}

\noindent and the numbers $V_{n}^{ \pm} ,V_{n\alpha} ^{ \pm}  $
are defined by the following recurrent formulae:

\begin{equation}
\label{eq8}
\alpha ^{2}V_{\alpha} ^{ \pm}  + \alpha \sum\limits_{n = 1}^{\alpha}
{V_{n\alpha} ^{ \pm}  + \sum\limits_{s = 1}^{\alpha - 1} {\left( {q_{\alpha
- s} V_{s}^{ \pm}  \pm p_{\alpha - s} \sum\limits_{n = 1}^{s} {V_{ns}^{ \pm
}} }  \right)} + q_{\alpha}  = 0}
\end{equation}

\begin{equation}
\label{eq9} \alpha (\alpha  - n)V_{n\alpha }^ \pm   +
\sum\limits_{s = n}^{\alpha  - 1} {(q_{\alpha  - s}  \mp n \cdot
p_{\alpha  - s} )V_{ns}^ \pm   = 0}
\end{equation}

\begin{equation}
\label{eq10}
\alpha V_{\alpha} ^{ \pm}  \pm \sum\limits_{s = 1}^{\alpha - 1} {V_{s}^{ \pm
} p_{\alpha - s} \pm p_{\alpha}  = 0}
\end{equation}

\noindent
and the sequence (2.4) admits double termwise differentiation. Then under
the conditions (\ref{eq7}) we get

\begin{equation}
\label{eq11}
f_{ \pm}  \left( {t,\mu}  \right) = \Psi ^{ \pm} \left( {t} \right)e^{ \pm
i\mu t} + \int\limits_{t}^{\infty}  {K^{ \pm} \left( {t,u} \right)e^{ \pm
i\mu u}du} \quad ,
\end{equation}

\noindent where$K^{ \pm} \left( {t,u} \right)$,$\Psi ^{ \pm}
\left( {t} \right)$ have a form
\begin{equation}
 \label{eq4}
 K^\pm  (t,u) = \frac{1}{{2i}}\sum\limits_{n = 1}^\infty
{\sum\limits_{\alpha  = n}^\infty  {V_{n\alpha }^ \pm  e^{ -
\alpha t}  \cdot e^{ - (u - t)n/2} ,\Psi ^ \pm  (t) = 1 +
\sum\limits_{n = 1}^\infty  {V_n^ \pm  e^{ - nt} } } } ,
\end{equation}

 So, it is proved the following

\underline {Lemma1:} The function $\Psi ^{ \pm} \left( {t}
\right)$and the kernel of the transformation operator of the
equation (\ref{eq6})$K\left( {t,u} \right),u \ge t$, attached to$
+ \infty $, with the potentials (\ref{eq7}) permits the
representation (2.9), in which the series

$ \sum\limits_{ = 1}^\infty  {n^2 \left| {V_n^ \pm  } \right|}
$,$\sum\limits_{ò = 1}^{\infty}
{\frac{{1}}{{n}}\sum\limits_{\alpha = n + 1}^{\infty} {\alpha
\left( {\alpha - n} \right)\left| {V_{n\alpha} ^{ \pm} } \right|}}
$,$\sum\limits_{ò = 1}^{\infty}  {n \cdot \left| {V_{nn}^{ \pm} }
\right|} $

\noindent  are convergent.

 Remark: In our case the kernel of the
operator of transformation $K^{ \pm }\left( {t,u} \right),u \ge
t$, at $ + \infty $, and function $ \Psi ^ \pm  (t) $ are
constructed effectively.

It is possible to get the equality [1]
\begin{equation}
\label{eq4}
 f_{n}^{ \pm}  \left( {t} \right) = S_{n}^{ \pm}  f_{ \mp}  \left( {t, \mp
i\frac{{n}}{{2}}} \right),
\end{equation}

\noindent
where
 $f_{n}^{ \pm}  \left( {t} \right) = \mathop {lim}\limits_{\mu \to \mp in/2}
\left( {in \pm 2\mu}  \right)f_{ \pm}  \left( {t,\mu}  \right)$.

Rewriting the equality (2.10) in the form

\begin{equation}
\label{eq12}
\sum\limits_{\alpha = n}^{\infty}  {V_{n\alpha} ^{ \pm} }  e^{ - \alpha t}
\cdot e^{nt/2} = S_{n}^{ \pm}  e^{ - nt/2}\left( {1 + \sum\limits_{m =
1}^{\infty}  {V_{m}^{ \mp}  e^{ - mt} + \sum\limits_{m = 1}^{\infty}
{\sum\limits_{\alpha = m}^{\infty}  {\frac{{V_{m\alpha} ^{ \mp} } }{{i\left(
{m + n} \right)}}e^{ - \alpha t}}} } }  \right)
\end{equation}

\noindent
and denoting by
\begin{equation}
\label{eq4} z^{ \pm} \left( {t + s} \right) = \sum\limits_{m =
1}^{\infty}  {S_{m}^{ \pm}  e^{ - \left( {t + s} \right)m/2}}
\end{equation}

 \noindent from (\ref{eq12}) we obtain
the Marchenko type equation

\begin{equation}
\label{eq13} K^ \pm  (t,s) = \Psi ^ \pm  (t)z^ \pm  (t + s) +
\int\limits_t^\infty  {K^ \mp  (t,u)z^ \pm  (u + s)du} .
\end{equation}

So, the following is proved

 \underline {Lemma2:} If the
coefficients $\overline {p} \left( {t} \right),\overline {q}
\left( {t} \right)$of the equation (\ref{eq6}) have the form
(\ref{eq7}), then at every$t \ge 0$, the kernel of the
transformation operator (2.9) satisfies to the equation of the
Marchenko type (\ref{eq13}) in which the transition function $z^{
\pm} \left( {t} \right)$ has the form (2.12) and the
numbers$S_{m}^{ \pm}  $ are defined by the equality (2.10), from
which it is obtained, that$ S_m^ \pm   = V_{mm}^ \pm $.

Note that from relation (\ref{eq10}) it is easy to obtain  useful
in further formulas $\Psi ^{ +} \left( {t} \right) \cdot \Psi ^{
-} \left( {t} \right) = 1$ and $\mathop {lim}\limits_{x \to
\infty} \Psi ^{ \pm} \left( {x} \right) = 1$. The coefficients are
reconstructed by the kernel of the transformation operator and
function with the help of the formulas

\begin{equation}
\label{eq14} \Psi ^{ \pm} \left( {t} \right) = J \pm
i\int\limits_{t}^{\infty} {\overline {p} \left( {u} \right)\Psi ^{
\pm} \left( {u} \right)du},
\end{equation}

\begin{equation}
\label{eq15} K^{ \pm} \left( {t,t} \right) = \pm
\frac{{1}}{{2}}\int\limits_{t}^{\infty} {\overline {q} \left( {u}
\right)\Psi ^{ \pm} \left( {u} \right)du} \mp i\overline {p}
\left( {t} \right)\Psi ^{ \pm} \left( {t} \right) \pm
i\int\limits_{ +} ^{\infty}  {\overline {p} \left( {u} \right)K^{
\pm }\left( {u,u} \right)du}.
\end{equation}

Hence the basic equation (\ref{eq13}) and the form of the
transition function (2.12) make natural the formulation of the
inverse problem of reconstruction potentials of the equation
(\ref{eq6}) by numbers$ S_n^ \pm $.

In this formulation, which employs the transformation operator method, an
important moment is the proof of uniqueness solvability of the basic
equation (\ref{eq13}).

\underline {Lemma3}. The homogenous equation

\begin{equation}
\label{eq16} g^{ \pm} \left( {s} \right) -
\int\limits_{0}^{\infty}  {z^{ \pm} \left( {u + s} \right)g^{ \mp}
\left( {u} \right)du} = 0
\end{equation}

\noindent
 has only trivial solution in the space $L_{2} \left( {R^{ +} } \right)$.

\noindent The proof of the lemma 3 is analogous to [5, p.198].

\underline {Lemma 4}: For each fixed $a,\left( {Ima \ge 0} \right)$ the
homogenous equation

\begin{equation}
\label{eq17} g^ \pm  (s) - \int\limits_t^\infty  {z^ \pm  (u + s -
2ai)g \mp (u)du}  = 0
\end{equation}

\noindent has only trivial solution in the space $L_{2} \left(
{R^{ +} } \right).$

 Proof: In the equation (\ref{eq2}) we
substitute $x$ by $x + a$, where $ {\mathop{\rm Im}\nolimits} a
\ge 0 $. Then we obtain the same equation with the coefficients
$p_{\alpha}  \left( {x} \right) = p_{\alpha}  \left( {x + a}
\right),$ $\,\,q_{\alpha}  \left( {x} \right) = q_{\alpha}  \left(
{x + a} \right)$ belonging to $Q_{ +} ^{2} $\textit{.} Let us
remark, that the function $e_{ \pm}  \left( {x + a,\lambda}
\right)$ are solutions of the equation

$ - y'' + 2\lambda p_a \left( x \right)y + q_a \left( x \right)y =
\lambda ^2 y $

\noindent that at $x \to \infty $ have the form

$e_{ \pm} \left( {x + a,\lambda}  \right) = e^{ \pm ia\lambda} e^{
\pm i\lambda x} + o\left( {1} \right)$.

Therefore the function $e_{ \pm} ^{a} \left( {x,\lambda}  \right) = e^{ \mp
ia\lambda} e_{ \pm}  \left( {x + a,\lambda}  \right)$ is also solution of
the type (\ref{eq3}).

Then let us denote by $\left\{ {\hat {S}_{n}^{ \pm}  \left( {a}
\right)} \right\}_{n = 1}^{\infty}  $ the spectral data of the
operator  $L \equiv - \frac{{d^{2}}}{{dx^{2}}} + 2\lambda p_{a}
\left( {x} \right) + q_{a} \left( {x} \right) - \lambda ^{2}$. \\
According to (1.4)

 $\begin{array}{l}
 \,\,\,\,\,\,\,\,\,\hat {S}_{n}^{ \mp}  \left( {a} \right)e_{ \pm} ^{a}
\left( {x, \pm {\raise0.7ex\hbox{${n}$} \!\mathord{\left/ {\vphantom {{n}
{2}}}\right.\kern-\nulldelimiterspace}\!\lower0.7ex\hbox{${2}$}}} \right) =
\mathop {lim}\limits_{\lambda \to \pm {\raise0.7ex\hbox{${n}$}
\!\mathord{\left/ {\vphantom {{n}
{2}}}\right.\kern-\nulldelimiterspace}\!\lower0.7ex\hbox{${2}$}}} \left( {n
\mp 2\lambda}  \right)e_{ \mp} ^{a} \left( {x,\lambda}  \right) = \\
 \,\,\,\,\, = \mathop {lim}\limits_{\lambda \to \pm {\raise0.7ex\hbox{${n}$}
\!\mathord{\left/ {\vphantom {{n}
{2}}}\right.\kern-\nulldelimiterspace}\!\lower0.7ex\hbox{${2}$}}} \left( {n
\mp 2\lambda}  \right)e^{ \pm ia\lambda} e_{ \mp}  \left( {x + a,\lambda}
\right) = e^{ia{\raise0.7ex\hbox{${n}$} \!\mathord{\left/ {\vphantom {{n}
{2}}}\right.\kern-\nulldelimiterspace}\!\lower0.7ex\hbox{${2}$}}}\hat
{S}_{n}^{ \mp}  e_{ \pm}  \left( {x + a, \pm {\raise0.7ex\hbox{${n}$}
\!\mathord{\left/ {\vphantom {{n}
{2}}}\right.\kern-\nulldelimiterspace}\!\lower0.7ex\hbox{${2}$}}} \right) =
\\
 = e^{ia{\raise0.7ex\hbox{${n}$} \!\mathord{\left/ {\vphantom {{n}
{2}}}\right.\kern-\nulldelimiterspace}\!\lower0.7ex\hbox{${2}$}}}\hat
{S}_{n}^{ \mp}  e^{ia{\raise0.7ex\hbox{${n}$} \!\mathord{\left/
{\vphantom {{n}
{2}}}\right.\kern-\nulldelimiterspace}\!\lower0.7ex\hbox{${2}$}}}e_{
\pm} ^{a} \left( {x, \pm {\raise0.7ex\hbox{${n}$}
\!\mathord{\left/ {\vphantom {{n}
{2}}}\right.\kern-\nulldelimiterspace}\!\lower0.7ex\hbox{${2}$}}}
\right) = \hat {S}_{n}^{ \mp}  e^{ian}e_{ \pm} ^{a} \left( {x, \pm
{\raise0.7ex\hbox{${n}$} \!\mathord{\left/ {\vphantom {{n}
{2}}}\right.\kern-\nulldelimiterspace}\!\lower0.7ex\hbox{${2}$}}}
\right)\,
 \end{array}
$

Hence   $\hat {S}_{n}^{ \mp}  \left( {a} \right) = \hat {S}_{n}^{
\mp}  e^{ian}. $

Now discussing as in the above, we obtain the basic equation of
the form (\ref{eq13}) with the transition function

 $
Z_a^ \pm  \left( t \right) = \sum\limits_{n = 1}^\infty  {S_n^ \pm
\left( a \right)e^{ - {\raise0.7ex\hbox{${nt}$} \!\mathord{\left/
 {\vphantom {{nt} 2}}\right.\kern-\nulldelimiterspace}
\!\lower0.7ex\hbox{$2$}}}  = } \sum\limits_{n = 1}^\infty  {S_n^
\pm  e^{ian} e^{ - {\raise0.7ex\hbox{${nt}$} \!\mathord{\left/
 {\vphantom {{nt} 2}}\right.\kern-\nulldelimiterspace}
\!\lower0.7ex\hbox{$2$}}}  = } Z^ \pm  \left( {t - 2ia} \right) $

\noindent From lemmas 3 and 4 follows

\underline {Theorem 3}: The potentials $\overline {p} \left( {t}
\right),\overline {q} \left( {t} \right)$\textit{} of the\textit{}
equation (\ref{eq6}), satisfying to the condition (\ref{eq7}) are
uniquely defined by the numbers$S_{n}^{ \pm}  $. \section {III.
Proof of Theorem 2.}

\underline {Necessity:} The necessity of the condition (1) has
been proved in the [1].

To prove the necessity of the condition (2) we firstly show, that
from the trivial solvability of the main equation (\ref{eq13}) at
$t = 0$ in the class of functions satisfying to inequality$\left\|
{g\left( {u} \right)} \right\| \le Ce^{ - \frac{{u}}{{2}}},u \ge
0$, follows trivial solvability in $ l_2 (R^ +  ) $ of the
infinite system of equations

\begin{equation}
\label{eq18} g_{n}^{ \pm}  - \sum\limits_{m = 1}^{\infty}
{\frac{{2S_{m}^{ \pm} } }{{m + n}}g_{m}^{ \mp}  = 0} ,
\end{equation}

\noindent
where$g_{n}^{ \pm}  \in l_{2} \left( {R} \right);\,\,\,\,S_{n}^{ \pm}  \in
l_{1} $.

Rewrite (\ref{eq18}) in the form
\begin{equation}
\label{eq4}
 g_{n}^{ \pm}  - \sum\limits_{m = 1}^{\infty}  {\sum\limits_{k = 1}^{\infty}
{\frac{{4S_{m}^{ \mp}  S_{k}^{ \pm} } }{{\left( {m + k}
\right)\left( {n + k} \right)}}g_{m}^{ \pm}  = 0}}.
 \end{equation}

 Really, if $\left\{ {g_{n}}  \right\}_{n = 1}^{\infty}
\in l_{2} $ is a solution of this system, then the function

\begin{equation}
\label{eq19} g^ \pm  (u) =  - \sum\limits_{m = 1}^\infty
{\sum\limits_{k = 1}^\infty  {\frac{{2S_m^ \mp  S_k^ \pm  }}{{m +
k}}e^{ - ku/2} g_m^ \pm  } }
\end{equation}

\noindent
is defined for all $u \ge 0$, satisfies inequality

 $\left| {g^{ \pm} \left( {u} \right)} \right| \le c \cdot e^{ -
u/2};\,\,\,\,\,u \ge 0
$ and is a solution of (\ref{eq16}), as

 $
\begin{array}{l}
 g^ \pm  (s) - \int\limits_0^\infty  {\int\limits_0^\infty  {z^ \pm  (u + s)z^ \mp  (u + s_1 )g^ \pm  (s_1 )ds_1 du} }  =  \\
  - \sum\limits_{m = 1}^\infty  {\sum\limits_{k = 1}^\infty  {\frac{{2S_m^ \mp  S_k^ \pm  }}{{m + k}}e^{ - ks/2} g_m^ \pm  } }  + \int\limits_0^\infty  {\int\limits_0^\infty  {\sum\limits_{m = 1}^\infty  {\sum\limits_{k = 1}^\infty  {S_m^ \mp  S_k^ \pm  e^{ - (u + s)k/2}  \cdot e^{ - (u + s_1 )m/2} } } } }  \times  \\
  \times \left( {\sum\limits_{n = 1}^\infty  {\sum\limits_{r = 1}^\infty  {\frac{{2S_n^ \mp  S_r^ \pm  }}{{n + r}}e^{ - rs_1 /2} g_n^ \pm  } } } \right)ds_1 du =  - \sum\limits_{m = 1}^\infty  {\sum\limits_{k = 1}^\infty  {\frac{{2S_m^ \mp  S_k^ \pm  }}{{m + k}}e^{ - ks/2} g_m^ \pm   + } }  \\
  + \sum\limits_{m = 1}^\infty  {\sum\limits_{k = 1}^\infty  {\frac{{2S_m^ \mp  S_k^ \pm  }}{{m + k}}e^{ - ks/2} (\sum\limits_{n = 1}^\infty  {\sum\limits_{r = 1}^\infty  {\frac{{4S_n^ \mp  S_r^ \pm  }}{{(n + r)(m + r)}}} } g_n^ \pm  ) = } }  \\
  =  - \sum\limits_{m = 1}^\infty  {\sum\limits_{k = 1}^\infty  {\frac{{2S_m^ \mp  S_k^ \pm  }}{{m + k}}e^{ - ks/2} \left[ {g_m^ \pm   - \sum\limits_{n = 1}^\infty  {\sum\limits_{r = 1}^\infty  {\frac{{2S_r^ \pm  S_n^ \mp  }}{{(m + r)(n + r)}}} } g_n^ \pm  } \right] = 0} }  \\
 \end{array}$

  Since,$g^{ \pm} \left( {u} \right) = 0$ therefore, $S_{m}^{ \mp}
S_{k}^{ \pm}  g_{m}^{ \pm}  = 0$ for all $m \ge 1,\,k \ge 1$, and
$g_{m}^{ \pm}  = 0,\,m \ge 1$ according to (3.2). Let us introduce
in the space $ l_2 $ operator$F_{2}^{ \pm}  \left( {t} \right)$,
given by matrix

\begin{equation}
\label{eq20} F_{mn}^{ \pm 2} \left( {t} \right) = \sum\limits_{k =
1}^{\infty} {\frac{{4S_{n}^{ \mp}  S_{k}^{ \pm} } }{{\left( {n +
k} \right)\left( {m + k} \right)}}e^{ - \left( {m + k}
\right)t/2}e^{ - \left( {n + k} \right)t/2},  n,m \in N}.
\end{equation}

Then, from $n^{2}S_{n}^{ \pm}  \in l_{1} $, we get
that$\sum\limits_{j,k = 1}^{\infty}  {\left| {\left( {F_{2}^{ \pm}
\varphi _{j} ,\varphi _{k}} \right)_{l_{2}} }  \right|} < \infty
$, i.e.$ F(t) $ is a kernel operator [10]. So there exists the
determinant $\Delta ^{ \pm} \left( {t} \right) = det\left( {E -
F_{2}^{ \pm } \left( {t} \right)} \right)$ of the operator $E -
F_{2}^{ \pm}  \left( {t} \right)$ related, as it is not difficult
to see, with the determinant $D^{ \pm} \left( {z} \right)$ from
the condition 2) of the theorem 2, by relation $ \Delta ^ \pm  ( -
iz) = \det (E - F_2^ \pm  ( - iz)) \equiv D^ \pm (z) $.

The determinant of the system (\ref{eq18}) is$D^{ \pm} \left( {0} \right)$, and the
determinant of analogous system with potential $p_{z} \left( {x} \right) =
p\left( {x + z} \right),\,q_{z} \left( {x} \right) = q\left( {x + z}
\right),\,\,Imz > 0$ is

 $D^ \pm  (z) = \det \left\| {\delta _{mn}  - \sum\limits_{k =
1}^\infty  {\frac{{4S_n^ \mp  (z)S_k^ \pm  (z)}}{{(m + k)(n +
k)}}} } \right\|_{m,n = 1}^\infty   = \left| {\delta _{mn}  -
\sum\limits_{k = 1}^\infty  {\frac{{4S_n^ \mp  S_k^ \pm  }}{{(m +
k)(n + k)}}} e^{i\frac{{m + k}}{2}z} e^{i\frac{{n + k}}{2}z} }
\right|_{m,n = 1}^\infty $.

\noindent Therefore to prove the necessity of the condition 2) of
Theorem 2 one should check that$\Delta ^{ \pm} \left( {0} \right)
= D^{ \pm} \left( {0} \right) \ne 0$. Really, the system
(\ref{eq18}) may be written in $l_{2} $ as

\noindent
 $g^{ \pm}  - F_{2}^{ \pm}  \left( {0} \right)g^{ \pm}  = 0$

 \noindent
As $F_{2}^{ \pm}  \left( {0} \right)$ is a kernel operator, the
Fredholm theory is applicable to it. In accordance with this
theory trivial solvability of the last equation is equivalent to
the fact that $ \det (E + F_2^ \pm  (0)) $ is not equal to zero
[10]. Necessity of the condition 2) is proved.

\underline {Sufficiency:} Let us study (\ref{eq13}) in detail. It
is known [5] that $K^{ \pm} \left( {t,s} \right)$ can be expressed
by $\Psi ^{ \pm} \left( {t} \right)$ and solutions $P^{ \pm}
\left( {t,s} \right)$,$Q^{ \pm} \left( {t,s} \right)$ of the
Marchenko type equations from (\ref{eq13}) by the replacement of $
\Psi ^ \pm  (t) $ by 1 and $ \pm i$.

Then

\begin{equation}
\label{eq21} K^{ \pm} \left( {t,s} \right) = \Psi ^{ \mp} \left(
{t} \right)\alpha ^{ \pm }\left( {t,s} \right) + \Psi ^{ \pm}
\left( {t} \right)\beta ^{ \mp} \left( {t,s} \right),
\end{equation}

\noindent
where

\begin{equation}
\label{eq22}
\alpha ^{ \pm} \left( {t,s} \right) = \frac{{1}}{{2}}\left[ {P^{ \pm} \left(
{t,s} \right) \mp iQ^{ \pm} \left( {t,s} \right)} \right]
\end{equation}

\begin{equation}
\label{eq23}
\beta ^{ \mp} \left( {t,s} \right) = \frac{{1}}{{2}}\left[ {P^{ \pm} \left(
{t,s} \right) \pm iQ^{ \pm} \left( {t,s} \right)} \right]
\end{equation}

\begin{equation}
\label{eq24} \left[ {\Psi ^ \pm  (t)} \right]^2  = \frac{{1 -
\int\limits_t^\infty  {\left[ {\alpha ^ \pm  (t,u) - \beta ^ \pm
(t,u)} \right]} du}}{{1 - \int\limits_t^\infty  {\left[ {\alpha ^
\mp  (t,u) - \beta ^ \mp  (t,u)} \right]du} }}
\end{equation}

\noindent from which we uniquivocally define $\Psi ^{ \pm} \left(
{t} \right)$. We also take into account that the sign of $\Psi ^{
\pm} \left( {t} \right)$ is fixed from condition $\mathop
{lim}\limits_{t \to \infty}  \Psi ^{ \pm }\left( {t} \right) = 1$.

Thus for further studies we should consider the following equations
\begin{equation}
\label{eq25} P^{ \pm} \left( {t,s} \right) = z^{ \pm} \left( {t +
s} \right) + \int\limits_{t}^{\infty}  {P^{ \mp} } \left( {t,u}
\right)z^{ \pm} \left( {u + s} \right)du
 \end{equation}

\begin{equation}
\label{eq25}
 Q^ \pm  (t,s) =  \pm iz^ \pm  (t + s) +
\int\limits_t^\infty  {Q^ \mp  } (t,u)z^ \pm  (u + s)du
\end{equation}

Rewrite (3.9) in the form

\begin{equation}
\label{eq26}
\begin{array}{l}
 P^{ \pm} \left( {t,s} \right) = z^{ \pm} \left( {t + s} \right) +
\int\limits_{t}^{\infty}  {z^{ \mp} } \left( {u + t} \right)z^{ \pm} \left(
{u + s} \right)du + \\
 \,\,\,\,\,\,\,\,\,\,\,\,\,\,\,\,\,\,\,\,\, + \int\limits_{t}^{\infty}
{\int\limits_{t}^{\infty}  {P^{ \pm} \left( {t,\tau}  \right)z^{ \mp} \left(
{u + \tau}  \right)z^{ \pm} \left( {u + s} \right)dud\tau} }  \\
 \end{array}
\end{equation}

Let us introduce in the space $l_{2} $ operator $F_{1}^{ \pm}  \left( {t}
\right)$ defined by the matrix

\begin{equation}
\label{eq27}
F_{mn}^{ \pm 1} = \frac{{2S_{n}^{ \pm} } }{{m + n}}e^{ - \left( {m + n}
\right)t/2};\,\,\,Ret > 0
\end{equation}

Multiplying the equation (\ref{eq26}) by $ e^{ - mu/2} $ and
integrating it over $s \in \left[ {t,\infty} \right)$ we obtain

\begin{equation}
\label{eq28}
p^{ \pm} \left( {t} \right) = F_{1}^{ \pm}  \left( {t} \right)e\left( {t}
\right) + F_{2}^{ \pm}  \left( {t} \right)e\left( {t} \right) + p^{ \pm
}\left( {t} \right)F_{2}^{ \pm}  \left( {t} \right),
\end{equation}

\noindent where $F_{1}^{ \pm}  \left( {t} \right),F_{2}^{ \pm}
\left( {t} \right)$ is defined by the matrix
(\ref{eq20}),(\ref{eq27}) and

 $e\left( {t} \right) = \left\{ {e^{ - nt/2}E} \right\}_{n = 1}^{\infty}
;\,\,\,p^{ \pm} \left( {t} \right) = \left\{
{\int\limits_{t}^{\infty}  {P^{ \pm} } \left( {t,u} \right)e^{ -
nu/2}du} \right\}_{n = 1}^{\infty}  ; $

\noindent As $ F_2^ \pm (t) $ is kernel operator for $t \ge 0$ and
the condition $ \det (E - F_2^ \pm  (t) \ne 0 $ is satisfied,
there exists bounded in $l_{2} $ inverse operator$R^{ \pm }\left(
{t} \right) = \left( {1 - F_{2}^{ \pm} \left( {t} \right)}
\right)^{ - 1}$.

As $ F_1^ \pm  (t)e(t),F_2^ \pm  (t)e(t) \in l_2 $, from
(\ref{eq28}) we get

\begin{equation}
\label{eq29}
p^{ \pm} \left( {t} \right) = R^{ \pm} \left( {t} \right)\left[ {F_{1}^{ \mp
} \left( {t} \right) + F_{2}^{ \pm}  \left( {t} \right)} \right]e\left( {t}
\right).
\end{equation}

Now, let's take$ < f,g > = \sum\limits_{n = 1}^{\infty}  {f_{n} g_{n}}  $.
Then (\ref{eq26}) gives

\begin{equation}
\label{eq30}
\begin{array}{l}
 P^ \pm  (t,s) =  < e(t),B^ \pm  (s) >  +  < e(t),A^ \pm  (s,t) >  +  < p^ \pm  (t),A^ \pm  (s,t) >  =  \\
  =  < e(t),B^ \pm  (s) >  +  < e(t),A^ \pm  (s,t) >  +  < R^ \pm  (t)(F_1^ \pm  (t) + F_2^ \pm  (t)e(t),A^ \pm  (s,t) >  =  \\
  =  < e(t),B^ \pm  (s) >  +  < (R^ \pm  (t)(F_1^ \pm  (t) + F_2^ \pm  (t)) + 1)e(t),A^ \pm  (s,t) >  \\
 \end{array}
\end{equation}

\noindent where $B^{ \pm} \left( {s} \right) = \left\{ {B_{m}^{
\pm}  \left( {s} \right) = S_{m}^{ \pm}  e^{ - ms/2},s > 0}
\right\}_{m = 1}^{\infty}  $ \\
and $ A^ \pm  (s,t) = \left\{ {A_m^
\pm  (s,t) = \sum\limits_{k = 1}^\infty  {\frac{{2S_m^ \mp  S_k^
\pm  }}{{m + k}}} e^{ - ks/2} \cdot e^{ - (m + k)t/2} ;s,t > 0}
\right\} $

Now suppose that the conditions of the theorem are satisfied. Let us define
the function $P^{ \pm} \left( {t,s} \right)$ by the equality (\ref{eq30}) at $0
\le t \le u$ according to given above considerations. Then for $u \ge t$ we
obtain

 $
\begin{array}{l}
 P^ \pm  (t,s) - \int\limits_t^\infty  {\int\limits_t^\infty  {P^ \pm  (t,\tau )z^ \mp  (u + \tau )z^ \pm  (u + s)dud\tau } }  =  \\
  < e(t),B^ \pm  (s) >  +  < (R^ \pm  (t)(F_1^ \pm  (t) + F_2^ \pm  (t)) + 1)e(t),A^ \pm  (s,t) >  -  \\
  - \int\limits_t^\infty  {\left( { < e(t),B^ \pm  (\tau ) >  +  < (R^ \pm  (t)(F_1^ \pm  (t) + F_2^ \pm  (t)) + 1)e(t),A^ \pm  (\tau ,t) > } \right)\left( { < e(\tau ),A^ \pm  (s,t) > } \right)d\tau  = }  \\
  =  < e(t),B^ \pm  (s) >  +  < e(t),A^ \pm  (s,t) >  +  < R^ \pm  (t)(F_1^ \pm  (t) + F_2^ \pm  (t)e(t),A^ \pm  (s,t) >  \\
  -  < F_1^ \pm  (t)e(t),A^ \pm  (s,t) >  -  < F_2^ \pm  (t)e(t),A^ \pm  (s,t) >  -  \\
  -  < R^ \pm  (t)F_2^ \pm  (t)(F_1^ \pm  (t) + F_2^ \pm  (t))e(t),A^ \pm  (s,t) >  =  \\
  =  < e(t),B^ \pm  (s) >  +  < e(t),A^ \pm  (s,t) >  = z^ \pm  (t + s) + \int\limits_ + ^\infty  {z^ \mp  (u + t)z^ \pm  (u + s)du}  \\
 \end{array} $

 For the  $
Q^ \pm  (t,s) $ we analogously obtain, that

 \noindent
 $Q^{ \pm} \left( {t,s}
\right) = \pm i < e\left( {t} \right),B^{ \pm} \left( {s} \right)
> \mp i < e\left( {t} \right),A^{ \pm} \left( {s,t} \right) > + <
Q^{ \pm} \left( {t} \right),A^{ \pm} \left( {s,t} \right) > $

\noindent
 where

 $Q^{ \pm} \left( {t} \right) = \pm iR^{ \pm} \left( {t} \right)\left[
{F_{1}^{ \pm}  \left( {t} \right) - F_{2}^{ \pm}  \left( {t} \right)}
\right]e\left( {t} \right),$ à $Q^{ \pm} \left( {t} \right) = \left\{
{\int\limits_{t}^{\infty}  {Q^{ \pm} \left( {t,u} \right)} e^{ - nu/2}du}
\right\}_{n = 1}^{\infty}  $.

Thus the following lemma was proved.

\underline {Lemma 5:} At any $t \ge 0$the kernel $ K^ \pm  (t,s) $
of the transformation operator and function $\Psi ^{ \pm} \left(
{t} \right)$ satisfies the main equation

\[
K^{ \pm} \left( {t,s} \right) = \Psi ^{ \pm} \left( {t} \right)z^{ \pm
}\left( {t + s} \right) + \int\limits_{t}^{\infty}  {K^{ \mp} \left( {t,u}
\right)z^{ \pm} \left( {u + s} \right)du} .
\]

The uniquely solvability of the main equation follows from Lemma 3. By
substitution it is not difficult to calculate, that the solution of the main
equation, indeed, is

 $K^{ \pm} \left( {t,u} \right) = \frac{{1}}{{2i}}\sum\limits_{n = 1}^{\infty
} {\sum\limits_{\alpha = n}^{\infty}  {V_{n\alpha} ^{ \pm}  e^{ -
\alpha t} \cdot e^{ - \left( {u - t} \right)n/2},\,\,\Psi ^{ \pm}
\left( {t} \right) = 1 + \sum\limits_{n = 1}^{\infty}  {V_{n}^{
\pm}  e^{ - nt}} \,}}  , \quad $ where the numbers $ V_{n\alpha }^
\pm  ,V_n^ \pm $ are defined by the numbers $S_{n}^{ \pm}  ,n \in
N$ from recurrent relations

 $V_{mm}^{ \pm}  = S_{m}^{ \pm}  ,\,\,\,V_{m,\alpha + m}^{ \pm}  = S_{m}^{ \pm
} \left( {V_{\alpha} ^{ \mp}  + \sum\limits_{n = = 1}^{\alpha}
{\frac{{V_{n\alpha} ^{ \mp} } }{{n + m}}}}  \right) $

\noindent Now we are in position to prove the main statement of
the theorem, namely that the coefficient $\overline {p} \left( {t}
\right)$and $ \bar q(t) $have the form (\ref{eq7}). First from the
formula (\ref{eq14}) and (\ref{eq15}) we find that the potentials
$\bar {p}\left( {t} \right),\bar {q}\left( {t} \right)$ have the
form

\[
\bar {p}\left( {t} \right) = \sum\limits_{n = 1}^{\infty}  {ip_{n}}  e^{ -
nt},\,\,\,\bar {q}\left( {t} \right) = \sum\limits_{n = 1}^{\infty}  { -
q_{n} e^{ - nt}} ,
\]

\noindent where the numbers $p_{n,} q_{n} $ are defined by the
relations (\ref{eq9})-(\ref{eq11}). After, in order to prove that
the numbers $p_{n,} q_{n} $, satisfy the condition $\sum\limits_{n
= 1}^\infty  {n \cdot \left| {p_n } \right| < \infty ,}
\sum\limits_{n = 1}^{\infty}  {\left| {q_{n}}  \right| < \infty} $
we demonstrate for the matrix elements $R_{mn} \left( {t}
\right)$, of the operator $R\left( {t} \right)$ correctness of the
estimations

 $\left\| {R_{mn}^{ \pm}  \left( {t} \right)} \right\| \le \delta _{mn} +
C_{0} \left| {S_{n}^{\ast} }  \right|, \quad \left\| {\frac{{d^j
}}{{dt^j }}R_{mn}^ \pm  (t)} \right\| \le C_j \left| {S_n^* }
\right|,$ where $S_{n}^{\ast}  = max\left( {S_{n}^{ -}  ,S_{n}^{
+} } \right)$ and $C_{j} ,j \in N$ are some constants. Really,
from the equality $ R^ \pm  (t) = E + R^ \pm  (t)F_2^ \pm  (t) $
follows that

\[
\begin{array}{l}
 \left\| {R_{mn}^ \pm  (t)} \right\| \le \delta _{mn}  + \left( {\sum\limits_{p = 1}^\infty  {\left\| {R_{mp}^ \pm  (t)} \right\|^2 } } \right)^{\frac{1}{2}} \left( {\sum\limits_{p = 1}^\infty  {\left\| {F_{pn}^{ \pm 2} (t)} \right\|^2 } } \right)^{\frac{1}{2}}  \le  \\
  \le \delta _{mn}  + 2 \cdot \left( {(R^ \pm  (t)R^{ \pm  * } (t)} \right)_{pp} \sum\limits_{p = 1}^\infty  {\left. {\left( {\sum\limits_{k = 1}^\infty  {\frac{{4S_n^ \mp  S_k^ \pm  }}{{(p + k)(n + k)}}} } \right)^2 } \right)^{{1 \mathord{\left/
 {\vphantom {1 2}} \right.
 \kern-\nulldelimiterspace} 2}}  \le }  \\
  \le \delta _{mn}  + 2 \cdot \left( {(R^ \pm  (t)R^{ \pm  * } (t)} \right)_{pp} \sum\limits_{p = 1}^\infty  {\left. {\frac{1}{{(p + 1)^2 }}\left( {\sum\limits_{k = 1}^\infty  {\left| {S_k^ \pm  } \right|} } \right)^2 } \right)^{{1 \mathord{\left/
 {\vphantom {1 2}} \right.
 \kern-\nulldelimiterspace} 2}} \left| {S_n^ \mp  } \right| \le }  \\
  \le \delta _{mn}  + \left\| {R(t)} \right\|_{l_{2 \to l_2 } } \left| {S_n^* } \right|. \\
 \end{array}
\]

From other hand, as it was noted, the operator-function $ R^ \pm
(t) $ exists and is bounded in $l_{2} $(as $F_{2}^{ \pm}  \left(
{t} \right)$is kernel operator by $t \ge 0$and$\Delta ^{ \pm}
\left( {t} \right) = det\left( {E + F_{2}^{ \pm}  \left( {t}
\right)} \right) \ne 0$) that proves the first inequality of
(\ref{eq24}).

To prove the estimation $ \frac{d}{{dt}}R^ \pm  (t) = R^ \pm
(t)F_2^ \pm  (t)R^ \pm  (t) $ (3.9) firstly we get the estimation

 $\quad
\begin{array}{l}
 \left\| {\frac{{d}}{{dt}}R_{mn}^{ \pm}  \left( {t} \right)} \right\| \le
\sum\limits_{q = 1}^{\infty}  {\sum\limits_{p = 1}^{\infty}  {\left\|
{R_{mq}^{ \pm}  \left( {t} \right)} \right\|\left\| {F_{qp}^{ \pm 2} \left(
{t} \right)} \right\|} \left\| {R_{pn}^{ \pm}  \left( {t} \right)} \right\|
\le \sum\limits_{q = 1}^{\infty}  {\sum\limits_{p = 1}^{\infty}  {\left(
{\delta _{mn} + C_{1} \left| {S_{q}^{\ast} }  \right|} \right) \times} } }
\\
 \times \left( {\left| {S_{p}^{\ast} }  \right|} \right)\left( {\delta _{pn}
+ C_{2} \left| {S_{n}^{\ast} }  \right|} \right) \le \left( {1 + C_{3}
\sum\limits_{p = 1}^{\infty}  {\left| {S_{p}^{\ast} }  \right|}}
\right)^{2} \cdot \left| {S_{n}^{\ast} }  \right| \le C_{4} \left|
{S_{n}^{\ast} }  \right| \\
 \end{array}$

 \noindent The estimation
 $\left\| {\frac{{d^{2}}}{{dt^{2}}}R_{mn}^{ \pm}  \left( {t} \right)}
\right\| \le C_{5} \left| {S_{n}^{\ast} }  \right|$
 can be proved analogously.

Using these estimations, from (\ref{eq13}) and (\ref{eq14}) one
can demonstrate the correctness of the estimations

 $\left\| {\frac{{d^2 }}{{dt^2
}}P^ \pm  (t,s)} \right\| \le \infty ,\left\| {\frac{{d^2 }}{{dt^2
}}Q^ \pm  (t,s)} \right\| \le \infty . $

Thus, the functions $K^{ \pm} \left( {t,s} \right)$and $\Psi ^{
\pm} \left( {t} \right)$ have the second derivatives over $t$.
From this we conclude that the series $\sum\limits_{n =
1}^{\infty}  {\alpha ^{2}} \left| {V_{\alpha} ^{ \pm} }  \right|$
and $ \sum\limits_{n = 1}^\infty  {n^{ - 1} \sum\limits_{\alpha  =
1}^n {(\alpha  + n)\left| {V_{n\alpha }^ \pm  } \right|} } $ are
convergence. The forms of the coefficients $\overline {p} \left(
{t} \right)$and $\overline {q} \left( {t} \right)$ are directly
determined from the form of the functions $K^{ \pm} \left( {t,s}
\right)$, $\Psi ^{ \pm} \left( {t} \right)$ employing the formulae
(1.15),(1.16) . We get that for the numbers $p_{n},q_{n}$ the
recurrent relations (\ref{eq9})-(\ref{eq11}) are correct and hence
the series$\sum\limits_{n = 1}^{\infty}  {n \cdot \left| {p_{n}}
\right| < \infty \,\,,} \sum\limits_{n = 1}^{\infty}  {\left|
{q_{n}} \right| < \infty}  $ are converges. Let, finally, $\left\{
{S_{n}^{ \pm} } \right\}_{n = 1}^{\infty} $ be a spectral data set
of the operator $L$ with the constructed coefficients $p(x),q(x)
\in Q_ + ^2 $. For completing of the proof it remains to show,
that $\left\{ {S_{n}^{ \pm} }  \right\}_{n = 1}^{\infty}
$coincides with the initial numbers$\left\{ {\hat {S}_{n}^{ \pm} }
\right\}_{n = 1}^{\infty} $. This follows from the
equality$S_{n}^{ \pm}  = V_{nn}^{ \pm}  = \hat {S}_{n}^{ \pm}  $.

The theorem is proved.

References:

\textbf{\textit{1.} }\textit{Efendiev R.F.} Inverse problem for
one class differential operator of the second \noindent order.
\textit{Reports National Acad. of Scien. of Azerbaijan ¹ 4-6,
15-20, 2001}.

\textbf{\textit{2.} }\textit{Efendiev R.F.} Spectral analysis of a
class of non-self-adjoint differential \noindent operator pencils
with a generalized function\textit{. Teoreticheskaya i}
 \textit{Matematicheskaya Fizika,Vol.145,‡'1.pp.102-107,October,2005(Russian)}
\textit{Theoretical and Mathematical
Physics,145(1):1457-1461,(2005),(Engish).}

\textbf{\textit{3.} }\textit{Melin A}. Operator methods for inverse
scattering on the real line.
\textit{ Comm.Partial Differential Equations, 1985.10:pp.677-766}

\textbf{\textit{4.} }\textit{ Marchenko V.A}. Sturm-Liouville operators and
applications, \textit{Birkhauser, Basel. 1986}

\textbf{\textit{5.} }\textit{Jaulent M. And. M. Jean}: The inverse s-Wave
Scattering Problem for a class of Potentials Depending on Energy.
\textit{Comm.Math.Phys. 28.1972, 177-220.}

\textbf{\textit{6.} }\textit{Deift P and Trubowitz D}. Inverse scattering on
the line. \textit{Comm.Pure Appl.Math. 1979,32, pp.121-251.}

\textbf{\textit{7.} }\textit{Aktosun T, Klaus M}. Inverse theory:
problem on the line. \textit{Chapter 2.2.4 in: Scattering (ed E.R.
Pike and P.C.Sabatier, Academic Press,} \textit{London, 2001,).
pp.770-785.}

\textbf{\textit{8.} }\textit{Gasymov M.G}. Spectral analysis of a
class non-self-adjoint operator of the second order.\textit{
Functional analysis and its appendix. (In Russian) 1980 V14.¹1.pp.
14-19. }

\textbf{\textit{9.} }\textit{Pastur L.A., Tkachenko V.A}. An
inverse problem for one class of one dimensional Schrödinger's
operators with complex periodic potentials. \textit{Functional
analysis and its appendix. (In Russian) 1990 V54.¹6.pp. 1252-1269}

\textbf{\textit{10.} }\textit{ Smirnov V.I.} Course of higher
mathematics. \textit{Ò.4. M.GITTL} .

\textbf{\textit{11.} }\textit{Gohberg I.T., Crein M.G}. Introduction to the
theory of linear nonselfadjoint
 operators. \textit{Nauka .1965,p.448}

\end{document}